\documentclass{elsart3-1}

\usepackage{graphicx}
\usepackage{epstopdf}
\usepackage{xcolor}
\usepackage{amssymb}
\usepackage{amsmath}
\usepackage{longtable}

\usepackage[english,francais]{babel}
\newtheorem{theorem}{Theorem}[section]
\newtheorem{conjecture}{Conjecture}[section]
\newtheorem{lemma}[theorem]{Lemma}
\newtheorem{e-proposition}[theorem]{Proposition}
\newtheorem{corollary}[theorem]{Corollary}
\newtheorem{e-definition}[theorem]{Definition\rm}

\newtheorem{example}{\it Example\/}
\def\og{\leavevmode\raise.3ex\hbox{$\scriptscriptstyle\langle\!\langle$~}}
\def\fg{\leavevmode\raise.3ex\hbox{~$\!\scriptscriptstyle\,\rangle\!\rangle$}}

\begin{document}
% place in the next line the header (rubrique) chosen for your article,
% if you know it (you can also have 2, format : Header1/Header2
Number theory
\centerline{}
%%%\begin{frontmatter}
%\frontmatter

\thanks[label1]{Partially supported by Serbian Ministry of Education and Science, Project 174032.
}
% \author{Name\thanksref{label2}}
\selectlanguage{english}
\title{The number of nonunimodular
roots of a reciprocal
polynomial}

% The [label1] can be suppressed if there is only one address for all authors

\selectlanguage{english}
\author{Dragan Stankov},
\ead{dstankov@rgf.bg.ac.rs}
%\author[authorlabel2]{Author Name2}
%\ead{author.name2@email.address2}

\address{
Katedra Matematike RGF-a,
Faculty of Mining and Geology,
University of Belgrade,
Belgrade, \DJ u\v sina 7,
Serbia}
%\address[authorlabel2]{Address2}

% If you know the dates of reception, and acceptation you can put them now;
%  idem the name of the person presenting the Note

%%%%%\medskip
%\begin{center}
%{\small Received 21 Novenber 2019; accepted after revision +++++\\
%Presented by Jean-Pierre Serre}
%\end{center}

%%%\begin{abstract}
\selectlanguage{english}
\section*{Abstract}
% Text of abstract in English

We introduce a sequence $P_{d}$ of monic reciprocal polynomials with integer coefficients having the central coefficients fixed as well as the peripheral coefficients. We prove that the ratio between number of
nonunimodular roots of $P_{d}$ and its degree $d$ has a limit $L$  when $d$ tends to infinity. We show that if the coefficients of a polynomial can be arbitrarily large in modulus then $L$ can be arbitrarily close to $0$. It seems reasonable to believe that if the %peripheral
coefficients are bounded then the analogue of Lehmer’s Conjecture is true: either $L=0$ or there exists a gap so that $L$ could not be arbitrarily close to $0$.  We present an algorithm for calculation the limit ratio and a numerical method for its approximation.
We estimated the limit ratio for a family of polynomials. We calculated the limit ratio of polynomials correlated to many bivariate polynomials having small Mahler measure.

\selectlanguage{english}
% main text
%\section{Introduction}
\section{Introduction}
The Mahler measure $M(P)$ of a polynomial $P(x)=a_dx^d+a_{d-1}x^{d-1}+\cdots +a_1x+a_0\in \mathbb{Z}[x]$ having $a_d \ne 0$ and zeros $\alpha_1,\alpha_2,\ldots,\alpha_d$ is defined as \[M(P(x)) := |a_d|\prod_{j=1}^{d}
\max(1, |\alpha_j|).\]
Let $I(P)$ denote the number of complex
zeros of $P(x)$ which are $<1$ in modulus, counted with multiplicities.
Let $U(P)$ denote the number of zeros of $P(x)$ which are $=1$ in modulus,
(again, counting with multiplicities). Such zeros %of the absolute value equal to one
are called unimodular.
Let $E(P)$ denote the number of complex
zeros of $P(x)$ which are $>1$ in modulus, counted with multiplicities.
Then it is obviously that $I(P)+U(P)+E(P)=d$.
Pisot number can be defined as a real algebraic integer greater than 1 having the minimal polynomial $P(x)$ of degree $d$ such that $I(P)=d-1$.
The minimal polynomial of a Pisot number is called Pisot polynomial.
Salem number %can be defined as
is a real algebraic integer $> 1$ having the minimal polynomial $P(x)$ of degree $d$ such that $U(P)=d-2$, $I(P)=1$. %\cite{Sta,Zai}.
The minimal polynomial of a Salem number is called Salem polynomial. It is well known that Salem polynomial is reciprocal.

We say that a polynomial of degree $d$ is reciprocal if
$P(x) = x^dP(1/x)$.
If moduli of coefficients are small then a reciprocal polynomial has many unimodular roots. A Littlewood polynomial is a polynomial all of whose coefficients are $1$ or $-1$. Mukunda \cite{Muk}
showed that every self-reciprocal Littlewood polynomial of odd
degree at least 3 has at least 3 zeros on the unit circle. Drungilas \cite{Dru} proved that every
self-reciprocal Littlewood polynomial of odd degree $n \ge 7$ has at least $5$ zeros on the unit
circle and every self-reciprocal Littlewood polynomial of even degree $n \ge 14$ has at least
$4$ unimodular zeros. In \cite{BCFJ} two types of very special Littlewood polynomials are considered: Littlewood polynomials with one sign change in the sequence of coefficients and Littlewood
polynomials with one negative coefficient. The numbers $U(P)$ and $I(P)$ of such Littlewood
polynomials $P$ are investigated. In \cite{BEFL} Borwein, Erd\'{e}lyi, Ferguson and Lockhart showed that there exists a cosine polynomial $\sum_{m=1}^{N}\cos(n_m\theta) $
with the $n_m$ integral and all different so that the number of its real zeros in $[0,2\pi)$ is
$O(N^{9/10}(\log N)^{1/5})$ (here the frequencies $n_m = n_m(N)$ may vary with $N$). However, there
are reasons to believe that a cosine polynomial $\sum_{m=1}^{N}\cos(n_m\theta) $ always has many zeros in
the period.

Clearly,
if $\alpha_j$, is a root of a reciprocal $P(x)$ then $1/\alpha_j$ is also a root of $P(x)$ so that $I(P)=E(P)$.
Let $C(P)=\frac{I(P)+E(P)}{2n}$ be the ratio between the number of nonunimodular zeros of $P$ and its degree. Actually, it is the probability that a randomly chosen zero is not unimodular, and $C(P)=\frac{E(P)}{n}$.

Here we will investigate a special sequence of polynomials.
Let $n$, $k$, $l$ $a_0,a_1,\ldots,a_k$, $b_0,b_1,\ldots,b_l$ be integers such that $2n>k\ge 0$, $l\ge 0$, and let $P_{2n+2l}(x)$ be a monic, reciprocal polynomial with integer coefficients
\begin{equation}\label{eq:SeqPoly}
P_{2n+2l}=x^{n+l}\left(\sum
_{j=0}^l b_j\left(x^{n+j}+\frac{1}{x^{n+j}}\right) +a_0+\sum
_{j=1}^k a_j\left(x^j+\frac{1}{x^j}\right)\right).
\end{equation}
We should remark that we have already studied in \cite{Sta} the special case of \eqref{eq:SeqPoly} for $l=0$, $b_0=1$. We are looking for the sequence $P_{2n+2l}$ %of monic reciprocal polynomials with integer coefficients,
such that the ratio $C(P_{2n+2l})$ has a limit when $n$ tends to $\infty$ and $0<\lim_{n\rightarrow\infty}C(P_{2n+2l})<1$. If $P_{2n+2l}$ is a sequence of Salem polynomials then this limit is trivially $0$. Such sequences are well known: Salem (see \cite{Sal} Theorem IV, p.30) found a simple way, which we present in the following example, to construct infinite sequences of Salem numbers (and Salem polynomials) from Pisot numbers.

\begin{example}\label{ExmLimit}
  If $x^l+b_{l-1}x^{l-1}+\cdots+b_0$ is a Pisot polynomial, $b_l=1$, $k=a_0=0$, then \eqref{eq:SeqPoly} is a sequence of Salem polynomials.
\end{example}

The definition of the Mahler measure could be extended to polynomials in several
variables. We recall Jensen's formula which states that
$\int_0^{1} \log |P(e^{2\pi i\theta})|\d \theta= \log |a_0| + \sum_{j=1}^d \log\max (|\alpha_j|, 1)$
Thus
$$M(P) = \exp \left\{ \int_0^{1} \log |P(e^{2\pi i\theta})|\d \theta\right\},$$
so M(P) is just the geometric mean of $|P(z)|$ on the torus $T$.
Hence a natural
candidate for $M(F)$ is
$$M(F) = \exp\left\{\int_0^{1} \d \theta_1\cdots\int_0^{1}\log |F(e^{2\pi i\theta_1},\ldots,e^{2\pi i\theta_r})|\d \theta_r \right\}.$$
%The smallest known Mahler measures in two variables are (see \cite{Boy})
%$$M((x + 1)y^2 + (x^2 + x + 1)y + x(x + 1)) = 1.25542\ldots$$ and
%$$M(y^2 + (x^2 + x + 1)y + x^2) = 1.28573\ldots .$$
Boyd and Mossinghoff in \cite{BM} listed in a table $48$ bivariate polynomials having small Mahler measure. Here we calculated the limit ratio of polynomials correlated to bivariate polynomials quadratic in $y$ and added them to the table. Flammang studied some other measures, defined for agebraic integers, in \cite{F1a1} \cite{F1a2}. 

We will use the following theorem of Erd\H{o}s and Tur\'{a}n to prove our Theorem \ref{ThmLimit}.
\begin{theorem}\label{thmErdTur} (Erd\H{o}s, Tur\'{a}n) Let $F(x)=\sum_{k=0}^d a_kx^k\in \mathbb{C}[x]$ with $a_da_0 \ne 0$, and let
$$N(F;\alpha,\beta)=\#\{\textrm{roots}\;\;r\in \mathbb{C}\;\textrm{of}\;F\;\textrm{with}\;\alpha \le \arg(r)\le \beta\}.$$
Then for all $ 0\le \alpha<\beta \le 2\pi$,
$$\left|\frac{N(F;\alpha,\beta)}{d}-\frac{\beta-\alpha}{2\pi}\right|\le\frac{16}{\sqrt{d}}\left[\log\left(\frac{|a_0|+\cdots+|a_d|}{\sqrt{|a_0a_d|}}\right)\right]^{1/2}.$$
\end{theorem}

\section{The Limit Ratio}

The following Theorem \ref{ThmLimit} is a generalisation of Theorem 2.1 which we proved in \cite{Sta}, more precisely Theorem 2.1 can be obtained from the following Theorem \ref{ThmLimit} if we take $l=0$.

\begin{theorem}\label{ThmLimit}
If $k>0$ is an integer and $b_j=b_l-j$, $j=0,1,\ldots,l$ then for all fixed integers $a_j$, $j=1, \ldots,k$ there is a limit $C(P_{2n+2l})$ when $n$ tends to infinity.
\end{theorem}

\begin{pf}
The theorem will be proved if we show that $1-C(P_{2n+2l})$  has a limit when $n$ tends to $\infty$. Since $1-C(P_{2n+2l})=\frac{U(P_{2n+2l})}{2n+2l}$ we have to count the unimodular roots of $P_{2n+2l}(x)$. The equation $P_{2n+2l}(x)=0$ is equivalent to

\begin{equation}\label{eq:BEA}
x^{n+l}\sum
_{j=0}^l b_j\left(x^{n+j}+\frac{1}{x^{n+j}}\right) =x^{n+l}\left(-a_0-\sum
_{j=1}^k a_j\left(x^j+\frac{1}{x^j}\right)\right).
\end{equation}

%Let $B_l(x)=b_0x^l + b_1x^{l-1}+\cdots+b_l$
Let $B(x)$ be the polynomial on the left side and let $A(x)$ be the polynomial on the right side of the previous equation.

%If we use the substitution $x=e^{it}$ in the equation $P_{2n}(x)=0$ we get
%$$e^{i(n+l)t}\left(\sum_{j=0}^l 2 b_j\cos (n+j)t +a_0+\sum_{j=1}^k 2 a_j\cos jt\right)=0$$
%Since $e^{int}\ne 0$ it follows that the equation is equivalent to
%\begin{equation}\label{eq:CES}
%\sum_{j=0}^l b_j\cos (n+j)t =-\frac{a_0}{2}-\sum_{j=1}^k a_j\cos jt.
%\end{equation}

Since $b_j=b_{l-j}$, $j=0,1,\ldots,l$ we have
\begin{align*}
  B(x) & =x^{n+l}\sum
_{j=0}^l b_j\left(x^{n+j}+\frac{1}{x^{n+j}}\right) \\
   & =\sum
_{j=0}^l b_j\left(x^{2n+l+j}+{x^{l-j}}\right)\\
   & =\sum
_{j=0}^l b_{j}x^{2n+l+j}+\sum
_{j=0}^l b_{l-j}{x^{l-j}}\\
   & =\sum
_{j=0}^l b_j\left(x^{2n+l+j}+{x^{j}}\right)\\
   & =\sum
_{j=0}^l b_j{x^{j}}\left(x^{2n+l}+1\right)  \\
   & =\left(x^{2n+l}+1\right)\sum
_{j=0}^l b_j{x^{j}} \\
   & =\left(x^{2n+l}+1\right)x^{l/2}\sum
_{j=0}^l b_j{x^{j-l/2}} .\\
\end{align*}

%$$x^{n+l}\sum
%_{j=0}^l b_j\left(x^{n+j}+\frac{1}{x^{n+j}}\right) .$$
%$$=\sum
%_{j=0}^l b_j\left(x^{2n+l+j}+{x^{l-j}}\right).$$
%$$=\sum
%_{j=0}^l b_j\left(x^{2n+l+j}+{x^{j}}\right).$$
%$$=\sum
%_{j=0}^l b_j{x^{j}}\left(x^{2n+l}+1\right) .$$
%$$=\left(x^{2n+l}+1\right)\sum
%_{j=0}^l b_j{x^{j}} .$$
%$$=\left(x^{2n+l}+1\right)x^{l/2}\sum
%_{j=0}^l b_j{x^{j-l/2}} .$$

Finally it follows that
\begin{equation}\label{eq:LEvenOdd}
B(x) =x^{n+l}\left(x^{n+l/2}+\frac{1}{x^{n+l/2}}\right)\sum
_{j=0}^l b_j{x^{j-l/2}} .
\end{equation}

Since we have to find unimodular roots we use the substitution $x=e^{it}$ in the equation \eqref{eq:BEA}. If $l$ is even then we have
\begin{equation}\label{eq:LEven}
B(e^{it})=e^{i(n+l)t}2\cos[(n+l/2)t]\left(\sum_{j=0}^{l/2-1} 2 b_j\cos [(l/2-j)t]+b_{l/2}\right).
\end{equation}

If $l$ is odd then it follows from \eqref{eq:LEvenOdd}

\begin{equation}\label{eq:LOdd}
B(e^{it})=e^{i(n+l)t}2\cos[(n+l/2)t]\left(\sum_{j=0}^{(l-1)/2} 2 b_j\cos [(l/2-j)t]\right).
\end{equation}

From the substitution $x=e^{it}$ it follows that $x$ is unimodular if and only if $t$ is real so that we have to count the real roots of $B(e^{it})=A(e^{it})$, ($t\in [0,2\pi))$.
We denote with $E(t)$ the function defined by terms enclosed within brackets of \eqref{eq:LEven} or of \eqref{eq:LOdd} i.e.
%let $E(t):=\sum_{j=0}^{l/2-1} 2 b_j\cos (l/2-j)t+b_{l/2}$ when $l$ is even and let $E(t):=\sum_{j=0}^{(l-1)/2} 2 b_j\cos (l/2-j)t$ when $l$ is odd.
\begin{equation}\label{eq:EnvelDef}
   E(t) :=
  \begin{cases}
    \sum_{j=0}^{l/2-1} 2 b_j\cos (l/2-j)t+b_{l/2} & \text{if $l$ is even,} \\
    \sum_{j=0}^{(l-1)/2} 2 b_j\cos (l/2-j)t & \text{if $l$ is odd.}
  \end{cases}
\end{equation}
If $t\in \mathbb{R}$ we can divide the equation $B(e^{it})=A(e^{it})$ with $2e^{i(n+l)t}\ne 0$ and obtain
$$\cos[(n+l/2)t]E(t)=-a_0/2-\sum_{j=1}^{k}a_j\cos jt$$
Let $\Gamma$ be the graph of $E(t)$, let $\Gamma_1$ be the graph of $f_1(t)=\cos (n+l/2)t\;E(t)$. Obviously for all $n$ $\Gamma_1$ is settled between graphs of $E(t)$ and $-E(t)$ and in certain points $\Gamma$ touches $\Gamma_1$. For that reason we call $E(t)$ the envelope of $f_1(t)$. %, the function on the left side of equation \eqref{eq:CES},
Let $\Gamma_2$ be the graph of
\begin{equation}\label{eq:f2Def}
  f_2(t):=-a_0/2-\sum_{j=1}^k a_j\cos jt.
\end{equation}
Then $U(P)$ is equal to the number of intersection points of $\Gamma_1$ and $\Gamma_2$.
These intersection points are obviously settled between curves $y=-|E(t)|$ and $y=|E(t)|$.
Graph $\Gamma_2$ of the continuous function $f_2$ and graph $\Gamma$ are fixed i.e. they do not depend on $n$, therefore there are $r$ subintervals $I_j$, such that $r$ is a finite integer, $I_j=[\alpha_{j},\beta_{j}]$, $0<\beta_{j-1}<\alpha_j<\beta_j<\alpha_{j+1}<2\pi$, such that if $t\in I_j$ then $|f_2(t)|\le|E(t)|$, where $\alpha_{j}$, $\beta_{j}$ are solutions of
\begin{equation}\label{eq:AlphaBeta}
  |E(t)|=|f_2(t)|.
\end{equation}

Using Theorem \ref{thmErdTur} of Erd\H{o}s and Tur\'{a}n we obtain
$$\left|\frac{N(P_{2n+2l};\alpha_{j},\beta_{j})}{2n+2l}-\frac{\beta_{j}-\alpha_{j}}{2\pi}\right|\le\frac{16}{\sqrt{2n+2l}}\left[\log\left(\frac{2\sum
_{j=0}^l |b_j|+|a_0|+2\sum
_{j=1}^k |a_j|}{\sqrt{|b_0b_l|}}\right)\right]^{1/2}.$$
If we introduce a constant $$D:=\left[\log\left(\frac{2\sum
_{j=0}^l |b_j|+|a_0|+2\sum
_{j=1}^k |a_j|}{\sqrt{|b_0b_l|}}\right)\right]^{1/2}$$ then it follows that

$$\frac{\beta_{j}-\alpha_{j}}{2\pi}-\frac{16}{\sqrt{2n+2l}}D\le\frac{N(P_{2n+2l};\alpha_{j},\beta_{j})}{2n+2l}\le\frac{\beta_{j}-\alpha_{j}}{2\pi}+\frac{16}{\sqrt{2n+2l}}D.$$
If we summarize the previous inequalities for $j=1,2,\ldots,r$ then we get
$$\sum_{j=1}^r\frac{\beta_{j}-\alpha_{j}}{2\pi}-r\frac{16}{\sqrt{2n+2l}}D\le\sum_{j=1}^r\frac{N(P_{2n+2l};\alpha_{j},\beta_{j})}{2n+2l}\le\sum_{j=1}^r\frac{\beta_{j}-\alpha_{j}}{2\pi}+r\frac{16 }{\sqrt{2n+2l}}D.$$
Finally we have to notice that
$\sum_{j=1}^r N(P_{2n+2l};\alpha_{j},\beta_{j})=U(P_{2n+2l})$ and
find the limit when $n$ tends to infinity. Using the squeeze theorem it follows that
$$\lim_{n\rightarrow \infty}\frac{U(P_{2n+2l})}{2n+2l}= \sum_{j=0}^r\frac{\beta_{j}-\alpha_{j}}{2\pi}$$
because $\lim_{n\rightarrow \infty}r\frac{16}{\sqrt{2n+2l}}D=0.$
\end{pf}
\qed

Theorem \ref{ThmLimit} enable us to introduce the following
\begin{e-definition}
  Let the limit of $C(P_{2n+2l})$ when $n$ tends to infinity be called limit ratio and denoted $LC(P_{2n+2l})$.
\end{e-definition}

%\begin{e-definition}
%A part of the graph of $f_1(t)=\cos (n+l/2)t$ such that $(k-1)\pi/(n+l/2) \le t \le k\pi/(n+l/2)$, $k\in \mathbb{Z}$ is $k$-th branch of $\cos (n+l/2)t$. The interval $[(k-1)\pi/(n+l/2),k\pi/(n+l/2)]$ is the domain of the $k$-th branch.
%\end{e-definition}
%We can define branch of $\cos (n+l/2)t$ as part of the graph of $f_1(t)=\cos (n+l/2)t$ where $f_1(t)\ge 0$ increase
%Each branch of $\cos (n+l/2)t$ obviously has exactly one intersection point with the $t$-axis. Using the same technique as in \cite{Sta} we can prove that if $n$ is large enough then each branch of $\cos (n+l/2)t$ also has exactly one intersection point with $\Gamma_2$.

It is well known that $S_1(x)=x^4-x^3-x^2-x+1$ is a Salem polynomial having two real roots: a Salem number $\gamma>1$, $1/\gamma$ and two complex unimodular roots $\theta, \bar{\theta}$. Let $S_m(x)=x^4+b_{1,m}x^3+b_{2,m}x^2+b_{3,m}x+1$ be the Salem polynomial of the Salem number $\gamma^m$ so that its coefficients should be $b_{0,m}=b_{4,m}=1, $
\begin{equation}\label{eq:B1}
b_{1,m}=b_{3,m}=-(\gamma^m+1/\gamma^m+\theta^m+\bar{\theta}^m),
\end{equation}
\begin{equation}\label{eq:B2}
b_{2,m}=2+\theta^m\gamma^m+\theta^m/\gamma^m+\bar{\theta}^m\gamma^m+\bar{\theta}^m/\gamma^m.
\end{equation}
%Let $l=4, b_{0,1}=b_{4,1}=1, b_{1,1}=b_{2,1}=b_{3,1}=-1, k=0, a_{0}=1$

\begin{example}\label{ExmLimit0}
Let $T_{2n+8,m}$ denote
$$T_{2n+8,m}(x)=x^{n+4}\left(\sum
_{j=0}^4 b_{j,m}\left(x^{n+j}+\frac{1}{x^{n+j}}\right) +2\right).$$
\end{example}

\begin{theorem}\label{ThmLimit0}
If we use the notation introduced in the previous example then $$\lim_{m\rightarrow \infty}LC(T_{2n+8,m}(x))=0.$$
\end{theorem}

\begin{pf}
In this example $l=4$ is even, $k=0$, $a_0=2$. We have to use \eqref{eq:EnvelDef} to calculate the envelope:
$E_m(t)=2\cos(2t)+2b_{1,m}\cos t+b_{2,m}$. We have to solve \eqref{eq:AlphaBeta} that is equivalent with $E_m(t)=1$ or $E_m(t)=-1$. Since $\cos 2t=2\cos^2 t-1$ the equations are quadratic in $\cos(t)$, so that, solving $E_m(t)=\pm 1$, we take the solutions in $[-1,1]$. From $E_m(t)=1$ we get $\cos \alpha_{m}=\frac{1}{4}\left(-b_{1,m}-\sqrt{b_{1,m}^2-4b_{2,m}+12}\right)$.
From $E_m(t)=-1$ we get $\cos \beta_{m}=\frac{1}{4}\left(-b_{1,m}-\sqrt{b_{1,m}^2-4b_{2,m}+4}\right)$. It remains to calculate
$$\lim_{m\rightarrow \infty}\left(\cos \beta_{m}-\cos \alpha_{m}\right)=\lim_{m\rightarrow \infty}\frac{2}{\sqrt{b_{1,m}^2-4b_{2,m}+12}+\sqrt{b_{1,m}^2-4b_{2,m}+4}}.$$
To show that the last limit is $0$ it is sufficient to show that $b_{1,m}^2-4b_{2,m}$ tends to $+\infty$ when $m\rightarrow \infty$.
Using \eqref{eq:B1} and \eqref{eq:B2}
$$b_{1,m}^2-4b_{2,m}=\left(\gamma^{2m}-2\gamma^m\theta^m-2\gamma^m\bar{\theta}^m\right)+
\left(1/\gamma^{2m}-2\bar{\theta}^m/\gamma^m-2\theta^m/\gamma^m+\theta^{2m}+\bar{\theta}^{2m}+2\theta^m\bar{\theta}^m-6\right).$$
The terms inside the first pair of parentheses are equal to
$$\gamma^{m}(\gamma^{m}-2\theta^m-2\bar{\theta}^m)\ge \gamma^{m}(\gamma^{m}-4)$$
so that they tend to $+\infty$ when $m\rightarrow \infty$.
Since all terms inside the second pair of parentheses are bounded or tend to zero it follows that $ b_{1,m}^2-4b_{2,m}$ tends to $+\infty$ when $m\rightarrow \infty$.

%$$1/g^2*m-2*zz^m/g^m-2*z^m/g^m+g^2*m+z^2*m+zz^2*m-2*g^m*z^m-2*g^m*zz^m+2*z^m*zz^m-6$$
\end{pf}

To determine the envelope in the following theorem we need the following lemmas which can be easily proved.

\begin{lemma}\label{lemmaCOdd}
$$\sin \frac{t}{2}\left(\sum_{j=1}^{m} 2 \cos jt+1\right)=\sin \frac{(2m+1)t}{2}$$
\end{lemma}

\begin{pf}

$$\sin \frac{t}{2}\left(\sum_{j=1}^{m} 2 \cos jt+1\right)
=\sin \frac{t}{2}\left(\sum_{j=0}^{m} 2\cos jt-1\right)$$

$$=\sin \frac{t}{2}\left(2\sum_{j=0}^{m} \cos jt\right)-\sin \frac{t}{2}$$

$$=2\cos \frac{mt}{2}\sin\frac{(m+1)t}{2}-\sin \frac{t}{2}$$

$$=\sin \frac{(2m+1)t}{2}+\sin \frac{t}{2}-\sin\frac{t}{2}$$

$$=\sin \frac{(2m+1)t}{2}$$
\end{pf}

\begin{lemma}\label{lemmaCEven}
$$\sin \frac{t}{2}\left(\sum_{j=1}^{m} 2 \cos \frac{(2j-1)t}{2}\right)=\sin mt$$
\end{lemma}

\begin{pf}

The formula is obviously true for $m=1$ because $2\sin \frac{t}{2}\cos\frac{t}{2}=\sin t$.
We suppose that the formula is true for $m=k$ i.e.
$$\sin \frac{t}{2}\left(\sum_{j=1}^{k} 2 \cos \frac{(2j-1)t}{2}\right)=\sin kt.$$
Using the product-to-sum formula it follows that the formula is true for $m=k+1$:
$$\sin \frac{t}{2}\left(\sum_{j=1}^{k+1} 2 \cos \frac{(2j-1)t}{2}\right)=\sin kt+2\sin \frac{t}{2}\cos \frac{(2k+1)t}{2}=\sin kt+\sin(k+1)t-\sin kt=\sin (k+1)t.$$
We conclude by mathematical induction that the formula holds for every natural number $m$.
\end{pf}

\begin{theorem}\label{ThmEnvelB1}
  If $b_0=b_1=\cdots=b_l=1$ is valid in \eqref{eq:SeqPoly} then
\begin{equation}\label{eq:Env}
  E(t)=\frac{\sin{\frac{(l+1)t}{2}}}{\sin \frac{t}{2}}.
\end{equation}
\end{theorem}

\begin{pf}

If $l$ is even then \eqref{eq:EnvelDef} gives
%$$B(e^{it})=e^{i(n+l)t}2\cos[(n+l/2)t]\left(\sum_{j=0}^{l/2-1} 2 \cos [(l/2-j)t]+1\right).$$
$$E(t)=\sum_{j=0}^{l/2-1} 2 \cos [(l/2-j)t]+1.$$
If we change the index of summation $J:=l/2-j$ and then reverse the order of summation we get
\begin{equation}\label{eq:LEven1}
%B(e^{it})=e^{i(n+l)t}2\cos[(n+l/2)t]\left(\sum_{J=1}^{l/2} 2 \cos Jt+1\right).
E(t)=\sum_{J=1}^{l/2} 2 \cos Jt+1.
\end{equation}
Finally using Lemma \ref{lemmaCOdd} it follows that
%$$B(e^{it})=e^{i(n+l)t}2\cos[(n+l/2)t]\left(\frac{\sin \frac{(l+1) t}{2}}{\sin \frac{t}{2}}\right).$$
$$ \frac{\sin \frac{(l+1) t}{2}}{\sin \frac{t}{2}}.$$
If $l$ is odd then \eqref{eq:EnvelDef} gives
%$$B(e^{it})=e^{i(n+l)t}2\cos[(n+l/2)t]\left(\sum_{j=0}^{(l-1)/2} 2 \cos [(l/2-j)t]\right).$$
$$E(t)=\sum_{j=0}^{(l-1)/2} 2 \cos [(l/2-j)t].$$
If we change the index of summation $J:=1/2+l/2-j$ and then reverse the order of summation we get
\begin{equation}\label{eq:LEven1}
%B(e^{it})=e^{i(n+l)t}2\cos[(n+l/2)t]\left(\sum_{J=1}^{(l+1)/2} 2 \cos [(J-1/2)t]\right).
E(t)=\sum_{J=1}^{(l+1)/2} 2 \cos [(J-1/2)t].
\end{equation}
Finally using Lemma \ref{lemmaCEven} we get
%$$B(e^{it})=e^{i(n+l)t}2\cos[(n+l/2)t]\left(\frac{\sin{\frac{(l+1)t}{2}}}{\sin \frac{t}{2}}\right).$$
$$E(t)=\frac{\sin{\frac{(l+1)t}{2}}}{\sin \frac{t}{2}}. $$
\end{pf}

In \cite{BM} Boyd and Mossinghoff introduced the following
\begin{e-definition}
 Let $\varphi_A(x)$ denote the polynomial $(x^A-1)/(x-1)$, and write
$$P_{A,B}(x,y)=x^{\max(A-B,0)}(\varphi_A(x)+\varphi_B(x)y+x^{B-A}\varphi_A(x)y^2).$$
\end{e-definition}

\begin{example}
Let $H_{2n+2l}(x)$ denote
$$H_{2n+2l}(x)=x^{n+l}\left(\sum
_{j=0}^l \left(x^{n+j}+\frac{1}{x^{n+j}}\right) +1\right).$$
\end{example}
We can show that
$$H_{2n+2l}(x)=P_{l+1,1}(x,x^{n+l})/x^{l}.$$
It is convenient to substitute $l=m-1$ in the previous example.
\begin{theorem}\label{ThmFig1}
If $m$ is an integer greater than 1 then
  $$\frac{2}{\pi(2m+1)}\frac{\sin\frac{(m-1)\pi}{2m}}{\sin\frac{\pi}{2m}}<LC(H_{2n+2m-2}(x))<\frac{2}{6m-\pi}\frac{\sin\frac{(m-1)\pi}{2m}}{\sin\frac{\pi}{2m}}$$
\end{theorem}

\begin{pf}
Since $b_0=b_1=\cdots=b_{m-1}$ in the previous example, we can use Theorem \ref{ThmEnvelB1} to determine the envelope:
$E_m(T)=\frac{\sin{\frac{mT}{2}}}{\sin \frac{T}{2}}. $ We have to solve \eqref{eq:AlphaBeta} that is equivalent with $|E_m(T)|=1/2$, $T\in [0,2\pi]$ because $k=0$, $a_0=1$. If we substitute $T=2t$ it follows that we have to solve $$2|\sin mt|= \sin t,\;\;t\in [0,\pi]$$
because we have to determine the sum of length of all intervals where $2|\sin mt|<|\sin t|$ on $[0,\pi]$.
\begin{figure}[t]
  \centering
  \includegraphics[width=0.9\textwidth, height=0.8\textheight]{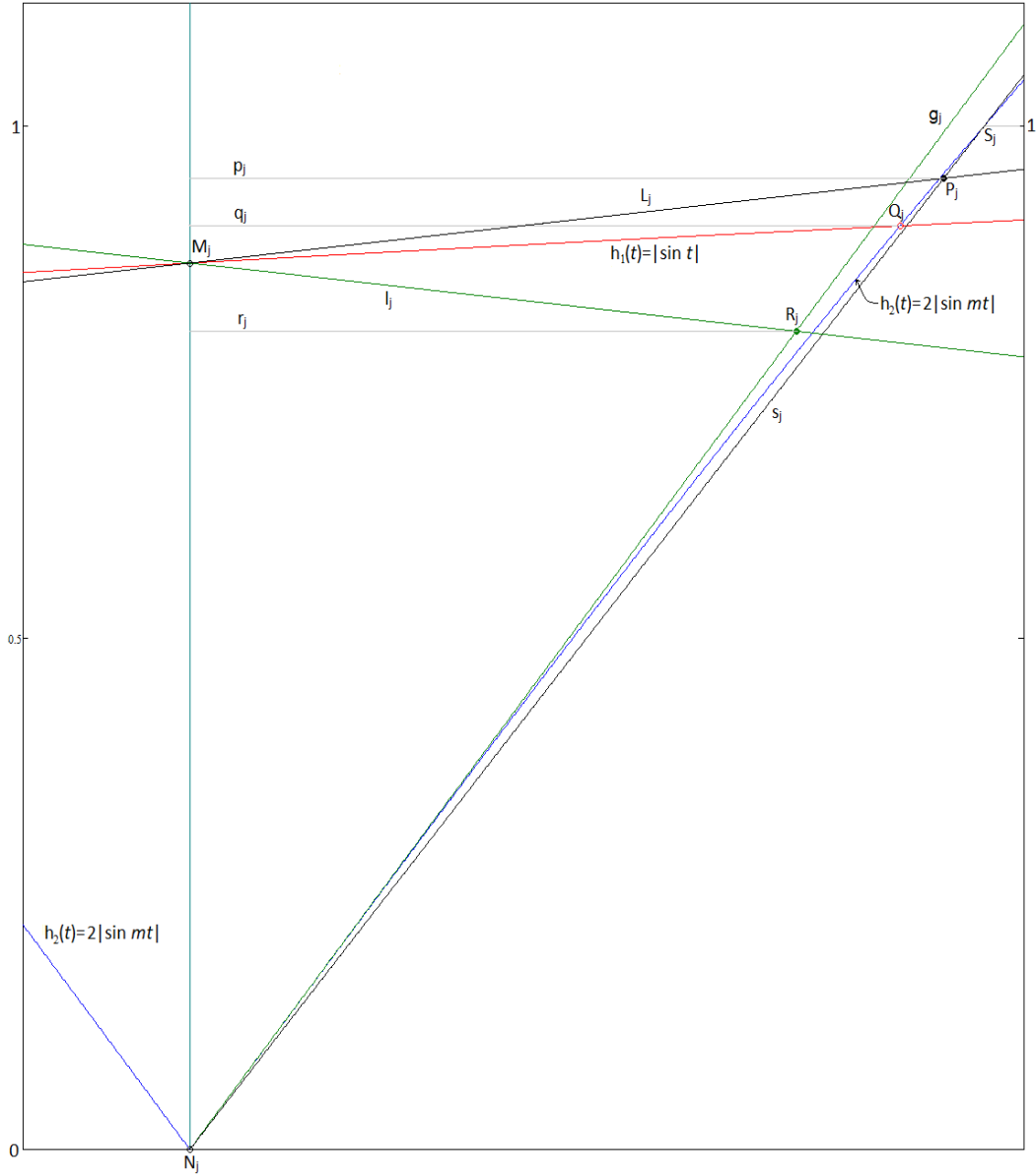}
\caption{
The estimation of the length of an interval where $2|\sin mt|<|\sin t|$.}\label{fig:Estim}       % Give a unique label
\end{figure}
Let $G_1$ be the graph of $h_1(t)=|\sin t|$ %, the function on the left side of equation \eqref{eq:CES},
and let $G_2$ be the graph of $h_2(t)=2|\sin mt|$. %The abscises of intersection points of $G_1$ and $G_2$ are equal to $\alpha_j/2$ or $\beta_j/2$, using the notation in the proof of the Theorem \ref{ThmLimit}. We are going to estimate $\beta_j-\alpha_j$.
Let $L_j$ be the line passing through $M_j\left(\frac{j\pi}{m},\sin\frac{j\pi}{m}\right)$ with the slope $1$, and let $l_j$ be the line passing through $M_j$ with the slope $-1$ (see. Fig. \ref{fig:Estim}). Let $g_j$ be the tangent line of $2|\sin mt|$ at $N_j\left(\frac{j\pi}{m},0\right)$ with the slope $2m$ and let $s_j$ be the secant line of $2|\sin mt|$ passing through $N_j$ and $S_j\left(\frac{j\pi}{m}+\frac{\pi}{6m},1\right)$. Let $Q_j$ be the unique intersection point of $G_1$ and $G_2$ on the segment $I_j=[\frac{j\pi}{m},\frac{j\pi}{m}+\frac{\pi}{2m}]$. Since $\frac{2}{\pi}<1$ there is the unique intersection point $P_j$ of $s_j$ and $L_j$, and also the unique intersection point $R_j$ of $g_j$ and $l_j$. On $I_j$ function $h_2$ increases and is concave down so that if $p_j$, $q_j$, $r_j$ are distances from points $P_j$, $Q_j$, $R_j$, respectively, to the vertical line $M_jN_j$ then $r_j<q_j<p_j$. To calculate $p_j$, $r_j$ it is convenient to use horizontal translation of all these objects such that $N_j$ moves to the origin $O$. Then $L_j$ moves to $L_j':y=t+\sin\frac{j\pi}{m}$ and $s_j$ moves to $s_j':y=\frac{6m}{\pi}t$. The solution $t$ of the system of these two equations is equal to $p_j$ so that
$$p_j=\frac{\sin\frac{j\pi}{m}}{\frac{6m}{\pi}-1}.$$
Also $l_j$ moves to $l_j':y=-t+\sin\frac{j\pi}{m}$ and $g_j$ moves to $g_j':y=2mt$. The solution $t$ of the system of last two equations is equal to $r_j$ so that
$$r_j=\frac{\sin\frac{j\pi}{m}}{2m+1}.$$

Similarly, let $\overline{Q}_j$ be the unique intersection point of $G_1$ and $G_2$ on the segment $\overline{I}_j=[\frac{j\pi}{m}-\frac{\pi}{2m},\frac{j\pi}{m}]$. Let $\overline{q}_j$ be the distances from point $\overline{Q}_j$ to the vertical line $M_jN_j$. Since the line $x=\frac{\pi}{2}$ is the axis of symmetry of $G_1$ as well as of $G_2$ it follows that $\overline{q}_j=q_{m-j}$ thus the sum of length of all intervals where $2|\sin mt|<|\sin t|$ on $[0,\pi]$ is equal to double $\sum_{j=1}^{m-1} q_j$. It follows from $r_j<q_j<p_j$ that $\frac{2}{\pi}\sum_{j=1}^{m-1} r_j<\frac{2}{\pi}\sum_{j=1}^{m-1} q_j<\frac{2}{\pi}\sum_{j=1}^{m-1} p_j$ so that
$$\frac{2}{\pi(2m+1)}\sum_{j=1}^{m-1} \sin\frac{j\pi}{m}<\frac{2}{\pi}\sum_{j=1}^{m-1} q_j<\frac{2}{6m-\pi}\sum_{j=1}^{m-1} \sin\frac{j\pi}{m}.$$
Finally if we use the formula for sum of sines with arguments in arithmetic progression we obtain the claim of the theorem.
\end{pf}
\qed

\begin{corollary}
If $A$ is an adherent point of the sequence $LC(P_{m,1}(x))_{m>1}$ then
  $$\frac{2}{\pi^2}\le A \le \frac{2}{3\pi}$$
\end{corollary}

\begin{pf}
We can easily show that the sequence of the lower bounds in the claim of previous theorem has the limit equal to $\frac{2}{\pi^2}\approx 0.2026$ and that the sequence of the upper bounds has the limit equal to $\frac{2}{3\pi}\approx 0.2122$ when $m\rightarrow \infty$.
\end{pf}

\begin{conjecture}
There is the limit of the sequence $LC(P_{m,1}(x))_{m>1}$ satisfying
\begin{equation}\label{eq:LimPm1}
\lim_{m\rightarrow \infty}LC(P_{m,1}(x))_{m>1}\approx 0.209.
\end{equation}
\end{conjecture}

\section{Approximating $\lim_{n\rightarrow\infty}C(P_{2n+2l})$}

It is necessary to explain how we approximated the limit in \eqref{eq:LimPm1}.
%\subsection{Algorithm for determination $\lim_{n\rightarrow\infty}C(P_{2n+2l})$}
In the proof of Theorem 1 we actually declared steps of an algorithm for determination $\lim_{n\rightarrow\infty}C(P_{2n+2l})$:

\begin{enumerate}
\item determine all real roots $t_j$ of the equations $f_2(t)=E(t)$ and $f_2(t)=-E(t)$, where $E(t)$, $f_2(t)$ are defined in \eqref{eq:EnvelDef} and \eqref{eq:f2Def},

\item arrange them as an increasing sequence $0=t_0<t_1<\ldots<t_p=2\pi$,

\item determine $r$ intervals $I_j=[\alpha_{j},\beta_{j}]$ such that if $\alpha_{j}<t<\beta_{j}$ then $|f_2(t)|\le |E(t)|$,  $\alpha_j,\beta_j \in \{t_0,t_1,\ldots,t_p\}$,

\item calculate $\lim_{n\rightarrow\infty}C(P_{2n+2l})=1-\sum_{j=1}^{r} (\beta_j-\alpha_j)/(2\pi).$
\end{enumerate}

If we bring to mind \eqref{eq:EnvelDef} it follows that the equation $f_2(t)=\pm E(t)$ i.e. $-a_0/2-\sum_{j=1}^k a_j\cos jt=\pm E(t)$ is algebraic in $\cos t$ so that $t_j$ can be expressed by arccosine of an algebraic real number $\alpha \in [-1,1]$ thus only solutions of this kind should be taken into account.

If $f_0(t)$ is defined:
\begin{equation*}
    f_0(t) = \begin{cases}
               1, & |f_2(t)|\ge |E(t)|\\
               0, & \text{otherwise}
           \end{cases}
\end{equation*}
then
\begin{equation}\label{eq:Cint}
\lim_{n\rightarrow\infty}C(P_{2n+2l})=\frac{1}{2\pi}\int_{0}^{2\pi}f_0(t)dt.
\end{equation}

We can approximate numerically the integral in \eqref{eq:Cint} i.e. $\lim_{n\rightarrow\infty}C(P_{2n+2l})$. Suppose the interval $[0,2\pi]$ is divided into $p$ equal subintervals of length $\Delta t=2\pi/p$ so that we introduce a partition of $[0,2\pi]$ $0=t_0<t_1<\ldots<t_p=2\pi$ such that $t_{j}-t_{j-1}=\Delta t$. Then we chose numbers $\xi_j\in [t_{j},t_{j-1}]$   and count all $\xi_j$ such that $|f_2(\xi_j)|\ge |E(t)$, $j=1,2,\ldots,p$. If there are $s$ such $\xi_j$ then $\lim_{n\rightarrow\infty}C(P_{2n+2l})$ is approximately equal to $\frac{s}{p}$.

$$\lim_{n\rightarrow\infty}C(P_{2n+2l})\approx\frac{1}{p}\sum_{j=1}^p f_0(j\frac{2\pi}{p})$$
where we chosed $\xi_j=2j\pi/p$.

If we introduce the substitution $t=2\pi u$ in \eqref{eq:Cint} we get

\begin{equation}\label{eq:Cint1}
\lim_{n\rightarrow\infty}C(P_{2n+2l})=\int_{0}^{1}f_0(2\pi u)du=\int_{U}du.
\end{equation}
where $U=\{u\in[0,1]:|f_2(2\pi u)|\ge |E(2\pi u)|\}$.

If we bring to mind the calculation of the Mahler measure in Exercise 2.24 and especially in Exercise 2.25 in the new book of McKee and Smyth \cite{MS}:
$$ M(P)=\exp \left( \int_{U}\log\frac{|f_2(2\pi u)|+\sqrt{f_2^2(2\pi u)-E^2(2\pi u)}}{|E(2\pi u)|}du \right)$$
where $U=\{u\in[0,1]:|f_2(2\pi u)|\ge |E(2\pi u)|\}$ then we can determine the correlation between Mahler measure and the limit ratio
$$LC(P)=\int_{U}du.$$

%\section{Limit points}

In Table \ref{tab:table0} we present $f_2(2\pi u)$ and $E(2\pi u)$ for certain families of polynomials, quadratic in $y$.
\begin{table}[h!]
  \begin{center}
    \caption{ $f_2(2\pi u)$ and $E(2\pi u)$ for certain families of polynomials.}
    \label{tab:table0}
    \begin{tabular}{l|c|l|l}
      Family & Definition & $f_2(2\pi u)$ & $E(2\pi u)$\\ % <-- added & and content for each column
      \hline
      $P_{a,b}(x,y)$ &  $x^{\max(a-b,0)}\left(\sum_{j=0}^{a-1}x^j+\sum_{j=0}^{b-1}x^jy+x^{b-a}\sum_{j=0}^{a-1}x^jy^2\right)$&  $\sin \left(\frac{b}{2}\pi u \right)$ & $2\sin \left(\frac{a}{2}\pi u \right)$\\ % <--
      $Q_{a,b}(x,y)$  &  $x^{\max(a-b,0)}(1 + x^a + (1 + x^b)y + x^{b-a}(1 + x^a)y^2$ &  $\cos \left(\frac{b}{2}\pi u\right)$& $2\cos \left(\frac{a}{2}\pi u\right)$ \\ % <--
      $R_{a,b}(x,y)$  &  $x^{\max(a-b,0)}(1 + x^a + (1 - x^b)y - x^{b-a}(1 + x^a)y^2$ &  $\sin \left(\frac{b}{2}\pi u\right)$& $2\cos \left(\frac{a}{2}\pi u \right)$\\ % <--
      $S_{a,b,\epsilon}(x,y)$  &  $1 + (x^a + \epsilon)(x^b + \epsilon)y + x^{a+b}y^2$, $\epsilon=\pm 1$ & $\cos \left(\frac{a+b}{2}\pi u\right)+\epsilon\cos\left(\frac{b-a}{2}\pi u\right)$ & 1\\ % <--
%      P &  &  & \\ % <--
    \end{tabular}
  \end{center}
\end{table}
In Table \ref{tab:table2} we present limit points calculated in \cite{BM} of Mahler measure of bivariate polynomials $P(x,y)$, quadratic in $y$, in ascending order. We complemented the table of Boyd and Mossinghoff by the limit points of the ratio between number of nonunimodular roots of the polynomial $P(x,x^n)$ and its degree when $n \rightarrow \infty$. As in \cite{BM} polynomials $P_{a,b}(x,y)$, $Q_{a,b}(x,y)$, $R_{a,b}(x,y)$, $S_{a,b,\epsilon}(x,y)$, defined in Table \ref{tab:table0}, are labeled as $P(a,b)$, $Q(a,b)$, $R(a,b)$, $S(a,b,\textrm{sgn}(\epsilon))$ respectively, in Table \ref{tab:table2}. Some polynomials are identified by the sequences, for example the third smallest known limit point $(1 + x) + (1 - x^2 + x^4)y + (x^3 + x^4)y^2$, is identified by [++000, +0$-$0+, 000++], as in \cite{BM}. Polynomials in Table \ref{tab:table2} are written explicitly in Table D.2 of \cite{MS}. We excluded the polynomials not quadratic in $y$.
It is interesting to compare Mahler measure and the limit ratio of polynomials in two variables.

\begin{enumerate}
  \item Mahler measure is $\ge 1$ while the limit ratio is in $[0,1]$.
  \item Mahler measures of two polynomials can be equal although its limit ratios are different (see examples (2) and (2') in Table \ref{tab:table2}.
  \item Mahler Measures of two polynomials increases although its limit ratio decreases.
  \item The polynomial $P_{2,3}$ has the smallest Mahler measure and the smallest limit ratio.
  \item The second smallest Mahler measure have $P_{2,1}$ and $P_{1,3}$ while the second smallest limit point has $R_{1,5}$.
\end{enumerate}

We showed in Example \ref{ExmLimit0} and Theorem \ref{ThmLimit0} that the limit ratio can be arbitrary close to zero. It is clear that in this example coefficients of the polynomials are unbounded. Our calculations show that if coefficients are bounded then the limit ratio can not be arbitrary close to zero. Also, the Theorem \ref{ThmFig1} supports our opinion that the analogue of Lehmer's conjecture is true:

\begin{conjecture}
  If $N$ is a natural number $\ge 1$ there is some $c(N)>0 $ such that any sequence $P_{2n+2l}$ of integer polynomials defined in \eqref{eq:SeqPoly}, having coefficients $\le N$ in modulus, that has the limit ratio strictly below $c(N)$ has the limit ratio equal to 0.
\end{conjecture}

  \begin{center}
    \begin{longtable}{l|c|l|l}
    \caption{ Limit points of Mahler measure and limit points of the ratio between number of nonunimodular roots of a polynomial and its degree.\label{tab:table2}}\\
%    \label{tab:table2}
%    \begin{tabular}{l|c|l|l}
      \textbf{$l$} Measure & Polynomial & $\lim\limits_{n\rightarrow\infty}C(S_l)$ & Exact value of $\lim\limits_{n\rightarrow\infty}C(S_l) $, sequence\\ % <-- added & and content for each column
      \hline
1. 1.2554338662666087457 & P(2, 3)& 0.1328095098966884 & $1-2\arccos(\frac{\sqrt{2}}{2}-\frac{1}{2})/\pi$\\
2. 1.2857348642919862749 & P(2, 1)& 0.1608612465103325 & $1-2\arccos(1/4)/\pi$\\
2'. 1.2857348642919862749 & P(1, 3)& 0.3333333333333333 & $1/3$\\
3. 1.3090983806523284595 &  & 0.2970136797597501 & [++000, +0$-$0+, 000++]\\
%4. 1.3156927029866410935 & P(3, 5)& 0.1646453474320021 & $\frac{4}{\pi}\arctan\left[\left(2\sqrt{94-26\sqrt{13}}+4\sqrt{13}-13\right)^{(-1)/2}\right]+$\\
%& & & $+\frac{-4}{\pi}\arctan\left[\left(\frac{1}{2}\sqrt{\frac{544\sqrt{5}}{121}+\frac{992}{121}}+\frac{4\sqrt{5}}{11}+\frac{17}{11}\right)^{(-1)/2}\right]$\\
4. 1.3156927029866410935 & P(3, 5)& 0.1646453474320021 & $\frac{4}{\pi}\arctan\frac{1}{\sqrt{2\sqrt{94-26\sqrt{13}}+4\sqrt{13}-13}}+$\\
& & & $+\frac{-4}{\pi}\arctan\frac{1}{\sqrt{\frac{1}{2}\sqrt{\frac{544\sqrt{5}}{121}+\frac{992}{121}}+\frac{4\sqrt{5}}{11}+\frac{17}{11}}}$\\
%5. 1.3247179572447460260 && 0 & T($1 + x - x^3$)\\
6. 1.3253724973075860349 & P(3, 4)& 0.1739784246485862 &\\%*
7. 1.3320511054374193142 & P(2, 5)& 0.2634504964561481 &\\%*
8. 1.3323961294587154121 & S(1, 3,+)& 0.3814904582918582 & $\arccos\left( \sqrt{17}/4-1/4\right)/\pi + 1/6$ \\
9. 1.3381374319388410775 & P(3, 2)& 0.1871346248477649 & $\frac{2}{\pi}\left[\pi-\arccos\frac{1+\sqrt{17}}{8}-\arccos\frac{1-\sqrt{17}}{8}\right] $\\
10. 1.3399999217381835332 & P(4, 7)& 0.1784746137157699 & \\
11. 1.3405068829308471079 & P(3, 1)& 0.1895159205822178 & $\frac{2}{\pi}\left[\arcsin(\sqrt{14}/4)-\arcsin(\sqrt{10}/4)\right]$\\
%12. 1.3497161046696958653 & & 0 & T($1 + x^2 - x^7$) [+++0$--$]\\
13. 1.3500148321630142650 & P(3, 7)& 0.2403097841316317 & \\
%14. 1.3503169790598690950 & S(1, 4,$-$)& 0.3105668890134219 & $\frac{2}{\pi}\left[\arccos(\beta)-\arccos(\alpha)\right]$, $\alpha,\beta$ \\
% & & & roots of $16z^5-24z^3+8z-1$ \\
15. 1.3511458956697046903 & P(4, 5)& 0.1902698620670582 &\\%*
16. 1.3524680625188602961 & P(5, 9)& 0.1860703555283188 & \\
17. 1.3536976494626355711 & Q(1, 6)& 0.1893226580984896 & \\
18. 1.3567481051456008311 & P(4, 3)& 0.1964065801899085 & \\
19. 1.3567859884526454967 & P(5, 8)& 0.1908351326172760 & \\
%20. 1.3581296324044179208 & & 0.37552129010217805 & [++, +0$-$$-$$-$0+,++]\\
20. 1.3581296324044179208 & & 0.3755212901021780 & $0.4-\alpha_1+0.8-2/3+1-\alpha_2$, $\alpha_1,\alpha_2$ \\
 & & & roots of $32z^6-48z^4+16z^2+2z-0.5$ \\
 & & & [++, +0$-$$-$$-$0+,++] \\
21. 1.3585455903960511404 & P(4, 1)& 0.1981783524823832 &\\%*
22. 1.3592080686995589268 & P(4, 9)& 0.2295536290347317 & \\
23. 1.3598117752819405021 & P(6, 11)& 0.1908185635976727 & \\
24. 1.3598158989877492950 & S(1, 6,+)& 0.3638326121576760 & \\
%25. 1.3599141493821189216 & & 0 & T($1 + x + x^8$) [+0$-$+0$-$+]\\
26. 1.3602208408592842371 & P(5, 7)& 0.1947758787175794 & \\
27. 1.3627242816569882815 & P(5, 6)& 0.1976969967166677 & \\
28. 1.3636514981864992177 & S(3, 5,+)& 0.3616163835316277 & \\
%29. 1.3641995455827723418 & & 0& T($1 - x^2 + x^5$)\\
%30. 1.3644358117806362770 & [+000, 00++, ++00, 000+]& & \\
31. 1.3645459857899151366 & P(7, 13)& 0.1940425569464528 & \\
32. 1.3646557293930641449 & P(5, 11)& 0.2236027778291241 & \\
33. 1.3650623157174417179 & S(2, 7,$-$)& 0.3360946113639976 & \\
34. 1.3654687370557201592 & P(5, 4)& 0.2007692138817449 & \\
%35. 1.3659850533667936783 & [++000, ++0−0, 00000, 0−0++, 000++]& & \\
36. 1.3661459663116649518 & P(5, 3)& 0.2014521139875612 & \\
37. 1.3665709746056369455 & P(5, 2)& 0.2018615118309531 & \\
38. 1.3668078899273126149 & P(5, 1)& 0.2020844014923849 & \\
39. 1.3668830708592258921 & R(1, 5)& 0.1417550822341309 & \\
40. 1.3669909125179202255 & P(7, 12)& 0.1970232013102869 & \\
41. 1.3677988580117157740 & P(8, 15)& 0.1963614081210482 & \\
%42. 1.3678546316653002345 && 0 & T($1 + x^4 + x^{11}$) [+$-$0+0$-$+0$-$+]\\
43. 1.3681962517212729703 & P(6, 13)& 0.2199360577499605 & \\
44. 1.3682140096679950123 & P(1, 9)& 0.2082012946810569 & \\
45. 1.3683434385467330804 & & 0.3045732337814742 & [++00000, ++0$-$0++,00000++]\\
46. 1.3687474425069274154 & P(6, 7)& 0.2014928273535877 & \\
47. 1.3689491694959833864 & P(7, 11)& 0.1994880038265199 & \\
48. 1.3697823199880122791 & S(1, 9,+)& 0.3622499773114010 & \\
%    \end{tabular}
\end{longtable}
  \end{center}

\end{document}